\documentclass[leqno]{article}
\usepackage{amssymb}
\usepackage{amsmath, enumerate, vmargin}

\makeatletter
\renewcommand\section{\@startsection {section}{1}{\z@}%
 {-3.5ex \@plus -1ex \@minus -.2ex}%
 {2.3ex \@plus.2ex}%
 {\center \normalfont\large\bfseries}}
\makeatother

\setmarginsrb{3cm}{3cm}{3cm}{3cm}{1mm}{6mm}{0mm}{10mm}

\newcommand{\Ec}{\mathbb{E}} 
\newcommand{\E}{\mathcal{E}} 
\newcommand{\M}{\mathcal{M}} 
\renewcommand{\H}{\mathcal{H}} 
 
\newcommand{\BMO}{\mathcal{BMO}} 
\newcommand{\MO}{\mathcal{MO}} 
\newcommand{\bmo}{\mathsf{bmo}}
\newcommand{\mo}{\mathsf{mo}} 
\newcommand{\h}{\mathsf{h}} 
 
\newcommand{\cqd}{\hfill$\Box$}

\newtheorem{theorem}{Theorem}[section]
\newtheorem{proposition}[theorem]{Proposition}

\newtheorem{remark}[theorem]{Remark}

\newenvironment{rk}{\begin{remark}\rm}{\end{remark}}

\newcommand{\pf}{\noindent{\it Proof.~~}}

\begin{document}

\title{A noncommutative Davis' decomposition for martingales}

\date{Mathilde Perrin}

\maketitle 

\begin{abstract}
We prove an analogue of the classical Davis' decomposition for martingales in noncommutative 
$L_p$-spaces, involving the square functions. 
We also determine the dual space of the noncommutative conditioned Hardy space $\h_1$. 
We further extend this latter result to the case $1<p<2$.
\end{abstract}

\makeatletter
 \renewcommand{\@makefntext}[1]{#1}
 
\makeatother \footnotetext{\noindent
M. Perrin: Laboratoire de Math{\'e}matiques, Universit{\'e} de France-Comt{\'e},
  25030 Besan\c{c}on Cedex,  France\\
  mathilde.perrin@univ-fcomte.fr\\
  Partially supported by the Agence Nationale de Recherche.
}
 
\section*{Introduction}

The theory of noncommutative martingale inequalities has been rapidly developed 
since the establishment of the noncommutative Burkholder-Gundy inequalities in \cite{px-BG}. 
Many of the classical martingale inequalities has been transferred to the noncommutative setting. 
These include, in particular, the Doob maximal inequality in \cite{ju-doob}, 
the Burkholder/Rosenthal inequality in \cite{jx-burk}, \cite{jx-ros}, 
several weak type $(1,1)$ inequalities in \cite{ran-mtrans,ran-weak,ran-cond} and the Gundy decomposition in \cite{par-ran-gu}. 
We would point out that the noncommutative Gundy's decomposition in this last work 
is remarkable and powerful in the sense that it implies several previous inequalities. 
For instance, it yields quite easily Randrianantoanina's weak type $(1,1)$ inequality on martingale transforms 
(see \cite{par-ran-gu}). 
It is, however, an open problem weather there exist a noncommutative analogue of the classical Davis' decomposition 
for martingales (see \cite{ran-cond} and \cite{par-burkh}). 
This is the main concern of our paper.

We now recall the classical Davis' decomposition for commutative martingales. 
Given a probability space $(\Omega, A, \mu)$, let $A_1, A_2, \cdots$ be an increasing filtration of $\sigma $-subalgebras of $A$ 
and let $\Ec_1,\Ec_2,\cdots$ denote the corresponding family of conditional expectations. 
Let $f=(f_n)_{n\geq 1}$ be a martingale adapted to this filtration and bounded in $L_1(\Omega)$. 
Then $M(f)=\sup|f_n|\in L_1(\Omega)$ iff 
we can decompose $f$ as a sum $f=g+h$ of two martingales adapted to the same filtration and satisfaying
$$s(g)=\Big(\displaystyle\sum_{n=1}^\infty\Ec_{n-1}|dg_n|^2\Big)^{1/2}\in L_1(\Omega) \quad 
\mbox { and } \quad \displaystyle\sum_{n=1}^\infty|dh_n| \in L_1(\Omega).$$
We refer to  \cite{gar} and \cite{dav} for more information. 

We denote by $\h_1$ the space of martingales $f$ with respect to $(A_n)_{n\geq 1}$ which admit such a decomposition 
and by $H_1^{\max}$ the space of martingales such that $M(f)\in L_1(\Omega)$. 
This decomposition appeared for the first time in \cite{dav} 
where Davis applied it to prove his famous theorem on the equivalence in $L_1$-norm between the martingale square function and Doob's maximal function:
$$\|M(f)\|_1 \approx \|S(f)\|_1$$
where $S(f)=\Big(\displaystyle\sum_{n\geq 1}|df_n|^2\Big)^{1/2}$. 
If we denote by $H_1$ the space of all $L_1$-martingales $f$ such that \\
$S(f) \in L_1(\Omega)$, 
then it turns out that the Hardy space $H_1$ coincides with the other two Hardy spaces:
$$H_1=\h_1=H_1^{\max}.$$

The main result of this paper is that the equality $H_1=\h_1$ holds in the noncommutative case. 
This answers positively a question asked in \cite{ran-cond}. 
This can be also considered as a noncommutative analogue of Davis' decomposition 
with the square function in place of the maximal function. 
Our approach to this result is via duality. 
We describe the dual space of $\h_1$ as a $\BMO$ type space. 
This is the second main result of the paper. 
Recall that this latter result is well known in the commutative case, 
the resulting dual of $\h_1$ is then the so-called small $\bmo$ (see \cite{prat}). 
Combining this duality with that between $\H_1$ and $\BMO$ established in \cite{px-BG}, 
we otain the announced equality $\H_1=\h_1$ in the noncommutative setting. 

Concerning $\H_1^{\max}$, it is shown in \cite{jx-const}, Corollary $14$, 
that $\H_1$ and $\H_1^{\max}$ do not coincide in general. 
More precisely $\H_1 \not\subset \H_1^{\max}$. 
But at the time of this writing we do not know if the reverse inclusion holds in the noncommutative setting.\\

The paper is organized as follows: 
in Section $1$ we give some preliminaries 
on noncommutative martingales and the noncommutative Hardy spaces. 
Section $2$ is devoted to the determination of the dual of $\h_1$, which allows us to show the equality $\H_1=\h_1$. 
This duality is extended to the case $1<p<2$ in Section $3$. 
There we describe the dual of $\h_p$ and use it to improve the estimation of an equivalence constant in the equivalence of the norms $\h_p$ and $\H_p$ given in \cite{ran-cond}.\\

After completing this paper, 
we learnt that Junge and Mei obtained the main result essentially at the same time 
(see Lemma $1.1$ of \cite{jm-riesz}). 
Note, however, that our proof of one direction in the duality theorem is different from theirs 
and yields a better constant (see Remark \ref{rk:duality-cst}).

\section{Preliminaries}

We use standard notation in operator algebras. We refer to \cite{kar-I} and \cite{tak-I} for background on von Neumann algebra theory. 
Throughout the paper all von Neumann algebras are assumed to be finite. 
Let $\M$ be a finite von Neumann algebra with a normal faithful normalized trace $\tau$. 
For $1\leq p \leq \infty$, we denote by $L_p(\M,\tau)$ or simply $L_p(\M)$ the noncommutative $L_p$-space associated with $(\M,\tau)$. 
Note that if $p=\infty$, $L_p(\M)$ is just $\M$ itself with the operator norm; also recall that the norm in $L_p(\M)$ ($1\leq p < \infty$) is defined as
$$\|x\|_p=(\tau(|x|^p))^{1/p}, \quad x\in L_p(\M)$$
where
$$|x|=(x^*x)^{1/2}$$
is the usual modulus of $x$. We refer to the survey \cite{px-survey} for more information on noncommutative $L_p$-spaces.\\

We now turn to the definition of noncommutative martingales. 
Let $(\M_n)_{n\geq 1}$ be an increasing sequence of von Neumann subalgebras of $\M$ such that the union of $\M_n$'s is weak$^*$-dense in $\M$. 
$(\M_n)_{n\geq 1}$ is called a filtration of $\M$. The restriction of $\tau$ to $\M_n$ is still denoted by $\tau$. 
Let $\E_n=\E(\,\cdot \,|\M_n)$ be the trace preserving conditional expectation of $\M$ with respect to $\M_n$. 
$\E_n$ defines a norm $1$ projection from $L_p(\M)$ onto $L_p(\M_n)$ 
for all $1\leq p \leq \infty$, and $\E_n(x)\geq 0$ whenever $x\geq 0$. 
A noncommutative martingale with respect to $(\M_n)_{n\geq 1}$ is a sequence $x=(x_n)_{n\geq 1}$ in $L_1(\M)$ such that 
$$\E_n(x_{n+1})=x_n, \quad \forall n\geq 1.$$
If additionally, $x\in L_p(\M)$ for some $1\leq p \leq\infty$, then $x$ is called an $L_p$-martingale. In this case, we set
$$\|x\|_p=\sup_{n\geq 1}\|x_n\|_p.$$
If $\|x\|_p<\infty$, then $x$ is called a bounded $L_p$-martingale. 
The difference sequence $dx=(dx_n)_{n\geq 1}$ of a martingale $x=(x_n)_{n\geq 1}$ is defined by
$$dx_n=x_n-x_{n-1}$$
with the usual convention that $x_0=0$. \\

We now describe Hardy spaces of noncommutative martingales. 
Following \cite{px-BG}, 
for $1\leq p < \infty$ and any finite sequence $a=(a_n)_{n\geq 1}$ in $L_p(\M)$, we set
$$\|a\|_{L_p(\M;\ell_2^c)}=\Big\|\Big(\sum_{n\geq 1}|a_n|^2\Big)^{1/2}\Big\|_p, \quad 
\|a\|_{L_p(\M;\ell_2^r)}=\Big\|\Big(\sum_{n\geq 1}|a_n^*|^2\Big)^{1/2}\Big\|_p.$$
Then $\|\cdot\|_{L_p(\M;\ell_2^c)}$ (resp. $\|\cdot\|_{L_p(\M;\ell_2^r)}$) defines a norm on the family of finite sequences of $L_p(\M)$. 
The corresponding completion is a Banach space, denoted by $L_p(\M;\ell_2^c)$ (resp. $L_p(\M;\ell_2^r)$). 
For $p=\infty$, we define $L_{\infty}(\M;\ell_2^c)$ (respectively $L_{\infty}(\M;\ell_2^r)$) 
as the Banach space of the sequences in $L_{\infty}(\M)$ such that $\displaystyle \sum_{n \geq 1} x_n^*x_n$ 
(respectively $\displaystyle \sum_{n \geq 1} x_nx_n^*$) converge for the weak operator topology. 
We recall the two square functions introduced in \cite{px-BG}. 
Let $x=(x_n)_{n\geq 1}$ be an $L_p$-martingale. We define 
$$S_{c,n}(x)= \Big(\displaystyle \sum_{k= 1}^n |dx_k|^2\Big)^{1/2} \quad \mbox{ and } \quad
S_{r,n}(x)= \Big(\displaystyle \sum_{k= 1}^n |dx_k^*|^2\Big)^{1/2}.$$
If $dx \in L_p(\M;\ell_2^c)$ (resp. $dx \in L_p(\M;\ell_2^r)$), we set
$$S_{c}(x)= \Big(\displaystyle \sum_{k\geq 1} |dx_k|^2\Big)^{1/2} \quad \Big(\mbox{ resp. } \quad
S_{r}(x)= \Big(\displaystyle \sum_{k\geq 1} |dx_k^*|^2\Big)^{1/2}\Big).$$
Then $S_c(x)$ and $S_r(x)$ are elements in $L_p(\M)$. 
Note that $dx \in L_p(\M;\ell_2^c)$ if and only if the sequence $(S_{c,n}(x))_{n\geq 1}$ is bounded in $L_p(\M)$. 
In this case
$$S_c(x)=\lim_{n\rightarrow \infty}S_{c,n}(x) \quad \mbox{ (relative to the weak}^*\mbox{-topology for } p=\infty).$$
The same remark applies to the row square function. \\

Let $1\leq p < \infty$. 
Define $\H_p^c(\M)$ (resp.  $\H_p^r(\M)$) to be 
the space of all $L_p$-martingales with respect to $(\M_n)_{n\geq 1}$ such that 
$dx \in L_p(\M;\ell_2^c)$ (resp. $dx \in L_p(\M;\ell_2^r)$), and set
$$\|x\|_{\mathcal{H}_p^c(\M)}=\|dx\|_{L_p(\M;\ell_2^c)} \quad \mbox{ and } \quad \|x\|_{\mathcal{H}_p^r(\M)}=\|dx\|_{L_p(\M;\ell_2^r)}.$$
Equipped respectively with the previous norms, $\mathcal{H}_p^c(\M)$ and $\mathcal{H}_p^r(\M)$ are Banach spaces.

Then we define the \emph{Hardy space of noncommutative martingales} as follows: \\
if $1\leq p <2$,
$$ \mathcal{H}_p(\M)=\mathcal{H}_p^c(\M)+\mathcal{H}_p^r(\M)$$
equipped with the sum norm
$$\|x\|_{\mathcal{H}_p(\M)}
=\inf \{ \|y\|_{\mathcal{H}_p^c(\M)}+\|z\|_{\mathcal{H}_p^r(\M)} : x=y+z, y \in \mathcal{H}_p^c(\M), z \in \mathcal{H}_p^r(\M)\};$$
if $2\leq p <\infty$,
$$ \mathcal{H}_p(\M)=\mathcal{H}_p^c(\M)\cap \mathcal{H}_p^r(\M)$$
equipped with the intersection norm
$$\|x\|_{\mathcal{H}_p(\M)}=\max\big(\|x\|_{\mathcal{H}_p^c(\M)}\; ,\; \|x\|_{\mathcal{H}_p^r(\M)}\big).$$

\vspace{0.5cm}

We now consider the conditioned versions of square functions and Hardy spaces developed in \cite{jx-burk}. 
Let $1\leq p <\infty$. For a finite $L_\infty$-martingale $x=(x_n)_{n\geq 1}$, define (with $\E_0=\E_1$)
$$\|x\|_{\h_p^c(\M)}=\Big\|\Big(\sum_{n=1}^{\infty} \E_{n-1}(|dx_n|^2)\Big)^{1/2} \Big\|_p$$
and
$$\|x\|_{\h_p^r(\M)}=\Big\|\Big(\sum_{n=1}^{\infty} \E_{n-1}(|dx_n^*|^2)\Big)^{1/2} \Big\|_p.$$
Let $\h_p^c(\M)$ and $\h_p^r(\M)$ be the corresponding completions. 
Then $\h_p^c(\M)$ and $\h_p^r(\M)$ are Banach spaces. 
We define the column and row conditioned square functions as follows. 
For any finite martingale $x=(x_n)_{n\geq 1}$ in $L_2(\M)$, set
$$s_{c}(x)= \Big(\displaystyle \sum_{n\geq1} \E_{n-1}(|dx_n|^2)\Big)^{1/2} \quad \mbox{ and } \quad
s_{r}(x)= \Big(\displaystyle \sum_{n\geq1} \E_{n-1}(|dx_n^*|^2)\Big)^{1/2}.$$
Then 
$$\|x\|_{\h_p^c(\M)}=\|s_c(x)\|_p \quad \mbox{and} \quad \|x\|_{\h_p^r(\M)}=\|s_r(x)\|_p.$$

We also need $\ell_p(L_p(\M))$, the space of all sequences $a=(a_n)_{n\geq 1}$ in $L_p(\M)$ such that
$$\|a\|_{\ell_p(L_p(\M))}=\Big(\sum_{n\geq 1}\|a_n\|_p^p\Big)^{1/p} <\infty.$$
Set
$$s_d(x)=\Big(\sum_{n\geq 1}|dx_n|_p^p\Big)^{1/p}.$$
We note that 
$$\|s_d(x)\|_p=\|dx\|_{\ell_p(L_p(\M))}.$$
Let $\h_p^d(\M)$ be the subspace of $\ell_p(L_p(\M))$ consisting of all martingale difference sequences.  

Following \cite{jx-burk}, we define the conditioned version of martingale Hardy spaces as follows:\\
if $1\leq p <2$,
$$ \h_p(\M)=\h_p^d(\M)+\h_p^c(\M)+\h_p^r(\M)$$
equipped with the norm
$$\|x\|_{\h_p(\M)}=\inf \{ \|x^d\|_{\h_p^d(\M)}+\|x^c\|_{\h_p^c(\M)}+\|x^r\|_{\h_p^r(\M)} \}$$
where the infimum is taken over all decompositions $x=x^d+x^c+x^r$ with $x^k \in \h_p^k(\M),k \in \{d,c,r\}$;\\
if $2\leq p <\infty$,
$$ \h_p(\M)=\h_p^d(\M)\cap \h_p^c(\M)\cap \h_p^r(\M)$$
equipped with the norm
$$\|x\|_{h_p(\M)}=\max\big(\|x\|_{\h_p^d(\M)} , \|x\|_{\h_p^c(\M)} , \|x\|_{\h_p^r(\M)}\big).$$

Throughout the rest of the paper letters like $\kappa_p, \nu_p \cdots$ 
will denote positive constants, which depend only on $p$ and may change from line to line. 
We will write $a_p\approx b_p$ as $p\rightarrow p_0$ to abbreviate the statement that there are two absolute positive constants $K_1$ and $K_2$ such that
$$K_1\leq \frac{a_p}{b_p}\leq K_2 \quad \mbox{ for } p \mbox{ close to }p_0.$$

\section{Noncommutative Davis' decomposition and the dual of $\h_1$}

Now we can state the main result of this section announced previously in introduction.

\begin{theorem}\label{th:H1=h1}
We have
$$\H_1(\M)=\h_1(\M) \quad \mbox{with equivalent norms.}$$
More precisely, if $x\in \H_1(\M)$,
$$\frac{1}{2}\|x\|_{\h_1}\leq \|x\|_{\H_1}\leq \sqrt{6}\|x\|_{\h_1}.$$
\end{theorem}

The inclusion $h_1(\M)\subset \H_1(\M)$ directly comes from the dual form of the reverse noncommutative Doob inequality in the case $0<p<1$ 
proved in \cite{jx-burk}, which is stated as follows. 
For all finite sequences $a=(a_n)_{n\geq 1}$ of positive elements in  $L_p(\M)$,
$$\Big\|\sum_{n\geq 1}a_n\Big\|_p\leq 2^{1/p} \Big\|\sum_{n\geq 1}\E_{n-1}(a_n)\Big\|_p.$$
Indeed, applying to $p=1/2$ and $a_n=|dx_n|^2$, we obtain for any martingale $x$ in $L_p$
$$\begin{array}{ccl}
\Big\|\Big(\displaystyle\sum_{n\geq 1}\E_{n-1}(|dx_n|^2)\Big)^{1/2}\Big\|_1&=&\Big\|\displaystyle\sum_{n\geq 1}\E_{n-1}(|dx_n|^2)\Big\|_{1/2}^{1/2}\\
&\geq &\displaystyle\frac{1}{4}\Big\|\displaystyle\sum_{n\geq 1}|dx_n|^2\Big\|_{1/2}^{1/2}\\
&=&\displaystyle\frac{1}{4}\Big\|\Big(\displaystyle\sum_{n\geq 1}|dx_n|^2\Big)^{1/2}\Big\|_1;
\end{array}$$
so $\|S_c(x)\|_1\leq 4 \|s_c(x)\|_1$. Similarly $\|S_r(x)\|_1\leq 4 \|s_r(x)\|_1$. 
On the other hand, we have
$$\begin{array}{ccl}
\|S_c(x)\|_1&=&\Big\|\displaystyle\sum_{n\geq 1}|dx_n|^2\Big\|_{1/2}^{1/2} \\
&\leq& \displaystyle\sum_{n\geq 1}\||dx_n|^2\|_{1/2}^{1/2}\\
&=&\displaystyle\sum_{n\geq 1}\|dx_n\|_1=\|s_d(x)\|_1.
\end{array}$$
Thus if $x\in h_1(\M)$, there exists $(x^d,x^c,x^r) \in h_1^d(\M)\times h_1^c(\M)\times h_1^r(\M)$ such that $x=x^d+x^c+x^r$ 
and from above $x^c,x^d \in \H_1^c(\M)$ and $x^r \in \H_1^r(\M)$, so $x\in \H_1(\M)$. 
Hence we deduce
$$\|x\|_{\H_1}\leq 4\|x\|_{\h_1}.$$

\vspace{0.5cm}

For the reverse inclusion, we will show the dual version. 
The dual approach gives also another proof for the direct inclusion, with a constant $\sqrt{6}$ instead of $4$. 
Recall that the dual space of $\H_1(\M)$ is the space $\BMO(\M)$ defined as follows (we refer to \cite{px-BG} for details). 
Set
$$\mathcal{BMO}^c(\M)=\{a \in L_2(\M): \sup_{n\geq1} \|\E_n(|a-\E_{n-1}(a)|^2)\|_{\infty} < \infty \}$$
where, as usual, $\E_{0}(a)=0$. $\BMO^c(\M)$ is equipped with the norm
$$\|a\|_{\mathcal{BMO}^c(\M)}=\Big(\sup_{n\geq1} \|\E_n(|a-\E_{n-1}(a)|^2)\|_{\infty}\Big)^{1/2}.$$
Then $(\mathcal{BMO}^c(\M),\|\cdot\|_{\mathcal{BMO}^c(\M)})$ is a Banach space. 
Similarly, we define 
$$\mathcal{BMO}^r(\M)=\{a \in L_2(\M): a^* \in \mathcal{BMO}^c(\M)\}$$
equipped with the norm
$$\|a\|_{\mathcal{BMO}^r(\M)}=\|a^*\|_{\mathcal{BMO}^c(\M)}.$$
Finally, we set
$$\mathcal{BMO}(\M)=\mathcal{BMO}^c(\M)\cap \mathcal{BMO}^r(\M)$$
equipped with the intersection norm
$$\|a\|_{\mathcal{BMO}(\M)}=\max \big(\|a\|_{\mathcal{BMO}^c(\M)},\|a\|_{\mathcal{BMO}^r(\M)}\big).$$
Note that if $a_n=\E_n(a)$, then 
$$\E_n(|a-\E_{n-1}(a)|^2)=\E_n\Big(\sum_{k\geq n} |da_k|^2\Big).$$

\vspace{0.5cm}

To describe the dual space of $\h_1(\M)$, we introduce similar spaces $\bmo^c(\M)$ and $\bmo^r(\M)$. 
Let
$$\bmo^c(\M)=\{a \in L_2(\M): \sup_{n\geq1} \|\E_n(|a-\E_{n}(a)|^2)\|_{\infty} < \infty \}$$
and equip $\bmo^c(\M)$ with the norm
$$\|a\|_{\bmo^c(\M)}=\max \Big(\|\E_1(a)\|_\infty \; , \; \Big(\sup_{n\geq1} \|\E_n(|a-\E_{n}(a)|^2)\|_{\infty}\Big)^{1/2}\Big).$$
This is a Banach space. 
Similarly, we define
$$\bmo^r(\M)=\{a \in L_2(\M): a^* \in \bmo^c(\M)\}$$
equipped with the norm
$$\|a\|_{\bmo^r(\M)}=\|a^*\|_{\bmo^c(\M)}.$$
For any sequence $a=(a_n)_{n\geq 1}$ in $\M$, we set
$$\|a\|_{\ell_\infty(L_\infty(\M))}=\sup_{n\geq 1}\|a_n\|_\infty.$$
Let $\bmo^d(\M)$ be the subspace of $\ell_\infty(L_\infty(\M))$ consisting of all martingale difference sequences.\\
Finally, we set
$$\bmo(\M)=\bmo^c(\M)\cap \bmo^r(\M)\cap \bmo^d(\M)$$
equipped with the intersection norm
$$\|a\|_{\bmo(\M)}=\max \big(\|a\|_{\bmo^c(\M)},\|a\|_{\bmo^r(\M)},\|a\|_{\bmo^d(\M)}\big).$$
Note that $\bmo^c(\M), \bmo^r(\M)$ and $\bmo(\M) \subset L_2(\M)$. 
As before, we have
$$\E_n(|a-\E_{n}(a)|^2)=\E_n\Big(\sum_{k>n}|da_k|^2\Big).$$
For convenience we denote $\mathcal{H}_1^c(\M), \mathcal{BMO}^c(\M), \h_1^c(\M),\bmo^c(\M) \cdots$, respectively, by $\mathcal{H}_1^c, \mathcal{BMO}^c,\h_1^c, \bmo^c \cdots
$\\
The relation between the spaces $\BMO$ and $\bmo$ can be stated as follows.

\begin{proposition}\label{prop:bmo}
We have 
$$\begin{array}{lcl}
\BMO^c&=&\bmo^c\cap \bmo^d,\\
\BMO^r&=&\bmo^r\cap \bmo^d,\\
\BMO&=&\bmo.
\end{array}$$
More precisely, for any $a\in L_2(\M)$,
$$ \|a\|_{\bmo^c\cap \bmo^d}\leq \|a\|_{\BMO^c} \leq \sqrt{2}\|a\|_{\bmo^c\cap \bmo^d}$$
and similar inequalities hold for the two other spaces.
\end{proposition}

\pf 
Let $a\in \BMO^c$. Then
$$\Big\|\E_n\Big(\sum_{k>n} |da_k|^2\Big)\Big\|_\infty \leq \Big\|\E_n\Big(\sum_{k\geq n} |da_k|^2\Big)\Big\|_\infty$$
and
$$\|da_n\|_\infty^2=\|\E_n|da_n|^2\|_\infty\leq \Big\|\E_n\Big(\sum_{k\geq n} |da_k|^2\Big)\Big\|_\infty.$$
Since $da_1=\E_1(a)$, taking the supremum over all $n\geq 1$ we find
$$\|a\|_{\bmo^c\cap \bmo^d}\leq\|a\|_{\BMO^c}.$$
Conversely, let $a\in \bmo^c \cap \bmo^d$, then
$$\Big\|\E_n\Big(\sum_{k\geq n} |da_k|^2\Big)\Big\|_\infty \leq 
\Big\|\E_n\Big(\sum_{k> n} |da_k|^2\Big)\Big\|_\infty+\|da_n\|_\infty^2.$$
Taking the supremum over all $n\geq 1$ we obtain
$$\|a\|_{\BMO^c}^2\leq \|a\|_{\bmo^c}^2+\|a\|_{\bmo^d}^2\leq 2\|a\|_{\bmo^c\cap \bmo^d}^2.$$
Hence
$$\|a\|_{\BMO^c} \leq \sqrt{2}\|a\|_{\bmo^c\cap \bmo^d}.$$
Passing to adjoints yields
$$\|a\|_{\bmo^r\cap \bmo^d}\leq \|a\|_{\BMO^r} \leq \sqrt{2}\|a\|_{\bmo^r\cap \bmo^d}.$$
These estimations show that the spaces $\BMO$ and $\bmo$ coincide.
\cqd

\vspace{0.5cm}

We have the following duality:

\begin{theorem}\label{th:dual_h1}
We have $(\h_1^c)^*=\bmo^c$ with equivalent norms. 
More precisely,
\begin{enumerate}[\rm (i)]
\item Every $a\in \bmo^c$ defines a continuous linear functional on $\h_1^c$ by
\begin{equation}\label{crochet dualite h1}
\phi_a(x)=\tau(a^*x), \quad \forall x\in L_2(\M).
\end{equation}
\item  Conversely, any $\phi \in (\h_1^c)^*$ is given as above by some $a \in \bmo^c$. \\
Moreover
$$\|a\|_{\bmo^c}\leq \| \phi_a\|_{(h_1^c)^*} \leq \sqrt{2}\|a\|_{\bmo^c}.$$
\end{enumerate}
Similarly, $(\h_1^r)^*=\bmo^r$ and $(\h_1)^*=\bmo$.
\end{theorem}

\vspace{0.5cm}

\begin{rk}\label{rk:dual_h1}
In the duality (\ref{crochet dualite h1}) we have identified an element $x \in L_2(\M)$ with the martingale $(\E_n(x))_{n\geq1}$. 
This martingale is in $\h_1^c$ and $\|x\|_{\h_1^c}\leq \|x\|_2$. 
Indeed, by the H\"older inequality, we have
$$\|x\|_{\h_1^c}=\|s_c(x)\|_1\leq \|s_c(x)\|_2=\|x\|_2,$$
where the last equality comes from the trace preserving property of conditional expectations 
and from the orthogonality in $L_2(\M)$ of martingale difference sequences. 
As finite $L_2$-martingales are dense in $\h_1^c$ and in $L_2(\M)$,  we deduce that $L_2(\M)$ is dense in $\h_1^c$.\\
\end{rk}

\pf 
\textbf{Step 1:}
We first show $\bmo^c\subset(h_1^c)^*$. 
This proof is similar to the corresponding one of the duality between  $\H_1$ and $\BMO$ in \cite{px-BG}. 
Let $a\in \bmo^c$. Define $\phi_a$ by (\ref{crochet dualite h1}). 
We must show that $\phi_a$ induces a continuous linear functional on $\h_1^c$.

Let $x$ be a finite $L_2$-martingale. 
Then (recalling our identification between a martingale and its limit value if the latter exists)
$$\phi_a(x)= \sum_{n\geq1}\tau(da_n^*dx_n).$$
Recall that 
$$s_{c,n}(x)=\Big(\sum_{k=1}^n\E_{k-1}|dx_k|^2\Big)^{1/2} 
\quad \mbox{ and } \quad 
s_{c}(x)=\Big(\sum_{k=1}^{\infty}\E_{k-1}|dx_k|^2\Big)^{1/2}.$$
By approximation we may assume that the $s_{c,n}(x)$'s are invertible elements in $\M$  for any $n\geq1$.\\
Then by the Cauchy-Schwarz inequality and the tracial property of $\tau$ we have
$$\begin{array}{ccl}
|\phi_a(x)|&=&\Big|\displaystyle \sum_{n\geq1} \tau( s_{c,n}(x)^{1/2}\,da_n^*\,dx_n\,s_{c,n}(x)^{-1/2})\Big| \\
&\leq&\Big[ \tau\Big(\displaystyle \sum_{n\geq1}s_{c,n}(x)^{1/2}|da_n|^2s_{c,n}(x)^{1/2}\Big)\Big]^{1/2}\\
&&\Big[ \tau\Big(\displaystyle \sum_{n\geq1}s_{c,n}(x)^{-1/2}|dx_n|^2s_{c,n}(x)^{-1/2}\Big)\Big]^{1/2}\\
&=& \Big[\tau\Big(\displaystyle \sum_{n\geq1}s_{c,n}(x)|da_n|^2\Big)\Big]^{1/2}
\Big[ \tau\Big(\displaystyle \sum_{n\geq1}s_{c,n}(x)^{-1}|dx_n|^2\Big)\Big]^{1/2} \\
&=&:I\cdot II.
\end{array}$$
To estimate $I$ we set $\theta_1=s_{c,1}(x)$ and $\theta_n=s_{c,n}(x)-s_{c,n-1}(x)$ for $n\geq 1$. 
Then $\theta_n \in L_1(\M_{n-1})$ and $s_{c,n}(x)=\displaystyle \sum_{k=1}^n \theta_k$. 
Using the Abel summation and the modular property of conditional expectations, we find
$$\begin{array}{ccl}
I^2&=&\displaystyle \sum_{n\geq1} \tau(s_{c,n}(x)|da_n|^2) 
= \displaystyle \sum_{n\geq1} \displaystyle \sum_{k=1}^n \tau(\theta_k|da_n|^2)\\
&=&\displaystyle \sum_{k\geq1}\tau \Big(\theta_k \displaystyle \sum_{n\geq k} |da_n|^2\Big) 
=\displaystyle \sum_{k\geq1}\tau \Big(\theta_k \E_{k-1}\Big(\displaystyle \sum_{n\geq k} |da_n|^2\Big)\Big) \\
&\leq&\displaystyle \sum_{k\geq1}\tau (\theta_k)\Big\|\E_{k-1}\Big(\displaystyle \sum_{n\geq k} |da_n|^2\Big)\Big\|_{\infty}\\
&\leq&\|x\|_{\h_1^c}\|a\|_{\bmo^c}^2.
\end{array}$$
To deal with $II$ first note that
$$ \tau\big[(s_{c,n}(x)^2-s_{c,n-1}(x)^2)s_{c,n}(x)^{-1}\big] 
= \tau\big[(s_{c,n}(x)-s_{c,n-1}(x))(1+s_{c,n-1}(x)s_{c,n}(x)^{-1})\big].$$
On the other hand, since $s_{c,n-1}(x)^2\leq s_{c,n}(x)^2$, we find
$$\begin{array}{ccl}
\|s_{c,n-1}(x)s_{c,n}(x)^{-1}\|_{\infty}^2&=&\|s_{c,n}(x)^{-1}s_{c,n-1}(x)^2s_{c,n}(x)^{-1}\|_{\infty}\\
&\leq &\|s_{c,n}(x)^{-1}s_{c,n}(x)^2s_{c,n}(x)^{-1}\|_{\infty}=1.
\end{array}$$
As $\E_{n-1}(|dx_n|^2)=s_{c,n}(x)^2-s_{c,n-1}(x)^2$ (with $s_{c,0}(x)=0$) we have
$$\begin{array}{ccl}
II^2&=&\displaystyle \sum_{n\geq1}\tau\big[\E_{n-1}(|dx_n|^2)s_{c,n}(x)^{-1}\big]\\
&=&\displaystyle \sum_{n\geq1} \tau\big[(s_{c,n}(x)-s_{c,n-1}(x))(1+s_{c,n-1}(x)s_{c,n}(x)^{-1})\big]\\
&\leq&\displaystyle \sum_{n\geq1} \tau\big[s_{c,n}(x)-s_{c,n-1}(x)\big] \|1+s_{c,n-1}(x)s_{c,n}(x)^{-1}\|_{\infty} \\
&\leq&2\tau\Big(\displaystyle \sum_{n\geq1} s_{c,n}(x)-s_{c,n-1}(x)\Big) \\
&=&2\tau(s_c(x))=2\|x\|_{\h_1^c}
\end{array}$$
Combining the preceding estimates on $I$ and $II$, we obtain, for any finite $L_2$-martingale $x$,
$$|\phi_a(x)| \leq \sqrt{2}\|x\|_{\h_1^c}\|a\|_{\bmo^c}.$$
Therefore $\phi_a$ extends to an element of $(\h_1^c)^*$ of norm $\leq \sqrt{2}\|a\|_{\bmo^c}$.\\
\smallskip\\
\textbf{Step 2:}
Let $\phi \in (\h_1^c)^*$ such that $\|\phi\|_{(\h_1^c)^*} \leq 1$. 
As $L_2(\M)\subset h_1^c$, $\phi$ induces a continuous functional $\tilde{\phi}$ on $L_2(\M)$. 
By the duality $(L_2(\M))^*=L_2(\M)$, there exists $a \in L_2(\M)$ such that 
$$ \tilde{\phi}(x)=\tau(a^*x), \quad \forall x \in L_2(\M).$$
By the density of $L_2(\M)$ in $h_1^c$ (see Remark \ref{rk:dual_h1}) we have 
\begin{equation}\label{norme phi}
\|\phi\|_{(\h_1^c)^*}=\sup_{x\in L_2(\M), \|x\|_{\h_1^c}\leq 1} |\tau(a^*x)|\leq 1.
\end{equation}
We will show that $a\in \bmo^c$. 
We want to estimate 
$$\|a\|_{\bmo^c}^2=\max\Big(\|\E_1(a)\|_\infty^2\; , \;\sup_n \|\E_n|a-\E_na|^2 \|_\infty\Big).$$
Let $x\in L_1(\M_1), \|x\|_1\leq 1$ be such that $\|\E_1(a)\|_\infty=|\tau(a^*x)|$. 
Then by \eqref{norme phi} we have 
$$\|\E_1(a)\|_\infty\leq  \|x\|_{\h_1^c}= \|x\|_1\leq 1.$$ 
On the other hand note that
$$\E_n|a-\E_na|^2 = \E_n\Big(\displaystyle\sum_{k>n}|da_k|^2\Big)
= \E_n\Big(\displaystyle\sum_{k>n}\E_{k-1}|da_k|^2\Big).$$
Fix $n\geq 1$. Let 
$$z=s_c(a)^2-s_{c,n}(a)^2=\sum_{k>n}\E_{k-1}|da_k|^2.$$
We note that $z\in L_1(\M)$ for $a\in L_2(\M)$ and the orthogonality of martingale difference sequences in $L_2(\M)$ gives
$$\|z\|_1=\tau(z)=\tau \Big(\sum_{k>n}|da_k|^2\Big)\leq \tau \Big(\sum_{k\geq 1}|da_k|^2\Big)=\tau(|a|^2) =\|a\|_2^2.$$
Let $x \in L_1^+(\M_n), \|x\|_1\leq 1$. 
Let $y$ be the martingale defined as follows 
$$dy_k=\left\{\begin{array}{ll}
0&\mbox{ if } k\leq n \\
da_kx&\mbox{ if } k> n 
\end{array}\right..$$
By (\ref{norme phi}) we have
$$\tau(a^*y) \leq \|y\|_{\h_1^c}.$$
Since $x\in L_1^+(\M_n)$, we have
$$\begin{array}{ccl}
\tau(a^*y)&=&\tau\Big(\displaystyle\sum_{k\geq 1} da_k^*dy_k\Big)
=\tau\Big(\displaystyle\sum_{k>n} |da_k|^2x\Big)\\
&=& \tau\Big(\displaystyle\sum_{k>n}\E_{k-1}(|da_k|^2x)\Big)\\
&=&\tau\Big(\displaystyle\sum_{k>n}\E_{k-1}(|da_k|^2)x\Big)=\tau(zx).
\end{array}$$
On the other hand, by the definition of $y$ and the fact that $x\in L_1^+(\M_n)$, we find
$$\begin{array}{ccl}
s_c(y)^2&=&\displaystyle\sum_{k\geq 1} \E_{k-1}|dy_k|^2
=\displaystyle\sum_{k>n} \E_{k-1}|da_kx|^2\\
&=&\displaystyle\sum_{k>n} \E_{k-1}(x|da_k|^2x)
=\displaystyle\sum_{k>n} x\E_{k-1}(|da_k|^2)x
=xzx.
\end{array}$$
Thus 
$$\|y\|_{\h_1^c}=\tau\Big((xzx)^{1/2}\Big).$$
Combining the preceding inequalities, we deduce
$$\tau(zx)\leq \tau((xzx)^{1/2}).$$
Since $x$ is positive, using the H\"{o}lder inequality, we find
$$\begin{array}{ccl}
\tau\Big((xzx)^{1/2}\Big) &=& \Big\|x^{1/2}(x^{1/2}zx^{1/2})x^{1/2}\Big\|_{1/2}^{1/2}\\
&\leq& \Big(\big\|x^{1/2}\big\|_2\big\|x^{1/2}zx^{1/2}\big\|_1\big\|x^{1/2}\big\|_2\Big)^{1/2}
= \tau(x)^{1/2}\tau(zx)^{1/2}.
\end{array}$$
It then follows that
$$\tau(zx)\leq \tau(x),$$
whence
\begin{equation}\label{ineg proj}
\tau(\E_n(z)x)=\tau(zx)\leq \tau(x)=\|x\|_1.
\end{equation}
Taking the supremum over all $x\in L_1^+(\M_n)$ with $\|x\|_1\leq 1$, we deduce $\|\E_n(z)\|_\infty\leq 1$. 
Therefore $a\in \bmo^c$ and $\|a\|_{\bmo^c}\leq 1$. 
This ends the proof of the duality $(\h_1^c)^*=\bmo^c$. 
Passing to adjoints yields the duality $(\h_1^r)^*=\bmo^r$.\\
\smallskip\\
\textbf{Step 3:}
Since finite martingales are dense in each $\h_1^c,\h_1^r$ and $\h_1^d$, 
the density property needed to apply the fact that the dual of a sum is the intersection of the duals holds. 
Thus it remains to determine the dual of $\h_1^d$. 
Since $\h_1^d$ is a subspace of $\ell_1(L_1(\M)$, the Hahn-Banach theorem gives
$$(\h_1^d)^*=\frac{(\ell_1(L_1(\M))^*}{(\h_1^d)^\perp}=\frac{\ell_\infty(L_\infty(\M))}{(\h_1^d)^\perp}.$$
Let 
$$P : \left\{
\begin{array}{ccc}
\ell_\infty(L_\infty)&\longrightarrow & \bmo^d\\
(a_n)_{n\geq1}&\longmapsto &(\E_n(a_n)-\E_{n-1}(a_n))_{n\geq1}
\end{array}\right..$$ 
We claim that $\ker P=(\h_1^d)^\perp$. 
Indeed, for $a\in \ker P$ and $x\in \h_1^d$ we have
$$\begin{array}{ccl}
\langle dx,a\rangle &=&\displaystyle\sum_{n\geq 1} \tau(dx_n^*a_n) 
=\displaystyle\sum_{n\geq 1} \Big[\tau(x_n^*a_n)-\tau(x_{n-1}^*a_n)\Big]\\
&=&\displaystyle\sum_{n\geq 1} \Big[\tau(x_n^*\E_n(a_n))-\tau(x_{n-1}^*\E_{n-1}(a_n))\Big] \\
&=& \displaystyle\sum_{n\geq 1} \tau(dx_n^*\E_{n-1}(a_n)) \quad \mbox{ for } \E_n(a_n)=\E_{n-1}(a_n)\\
&=&0.
\end{array}$$
Conversely, if $a \in (\h_1^d)^\perp$ 
we fix $n \geq 1$ and define the martingale $x$ by $dx_n=\E_n(a_n)-\E_{n-1}(a_n)$ and $dx_m=0$ if $m\neq n$. 
Since $a_n\in L_\infty(\M)$ and $\tau$ is finite, 
$$\sum_{m\geq 1}\|dx_m\|_1=\|dx_n\|_1 \leq 2\|a_n\|_1\leq 2\|a_n\|_\infty < \infty,$$
so $x \in \h_1^d$. Hence
$$\begin{array}{ccl}
0&=&\langle dx,a\rangle= \tau\big[(\E_n(a_n)-\E_{n-1}(a_n))^*a_n\big]\\
 &=&\tau(\E_n(a_n)^*\E_n(a_n))-\tau(\E_{n-1}(a_n)^*\E_{n-1}(a_n))\\
 &=&\tau |\E_n a_n-\E_{n-1} a_n|^2;
\end{array}$$
whence $\E_n(a_n)=\E_{n-1}(a_n)$. 
Thus we deduce that $a \in \ker P$. Therefore, our claim is proved. 
It then follows that $(\h_1^d)^*=\bmo^d$. 
Hence, the proof of the theorem is complete.
\cqd

\vspace{0.5cm}

We can now prove the reverse inclusion of Theorem \ref{th:H1=h1}.

\vspace{0.5cm}

\noindent{\it Proof of Theorem \ref{th:H1=h1}.~~} 
By the discussion following Theorem \ref{th:H1=h1}, we already know $\h_1\subset \H_1$. 
To prove the reverse inclusion, we use duality. 
It then suffices to show $(\h_1)^*\subset (\H_1)^*$. 
To this end, by Theorem \ref{th:dual_h1} and the duality theorem of \cite{px-BG}, we must show $\bmo \subset \BMO$. 
This result is stated in Proposition \ref{prop:bmo} , with the equivalence constant $\sqrt{2}$. 
Combining the estimation of Theorem \ref{th:dual_h1} and Proposition \ref{prop:bmo} with the appendix of \cite{px-BG}, 
we obtain for any $a\in (\h_1)^*$
$$\|a\|_{(\H_1)^*}\leq \sqrt{2}\,\|a\|_{\BMO} \leq  2 \|a\|_{\bmo} \leq 2\|a\|_{(\h_1)^*}$$
and
$$\|a\|_{(\h_1)^*}\leq \sqrt{2} \|a\|_{\bmo} \leq \sqrt{2} \|a\|_{\BMO} \leq \sqrt{6}\|a\|_{(\H_1)^*}.$$
\cqd

\vspace{0.5cm}

\begin{rk}
Combining Proposition \ref{prop:bmo} and the duality results, we also obtain
$$\H_1^c=\h_1^c+ \h_1^d \quad \mbox{ and } \quad \H_1^r=\h_1^r+ \h_1^d.$$ 
\end{rk}

\section{A description of the dual of $h_p$ for $1<p<2$}

In this section we extend the duality theorem in the previous section to the case $1<p<2$. 
Namely, we will describe the dual of $\h_p$ for $1<p<2$. 
The arguments are similar to those for $p=1$. 
The situation becomes, however, a little more complicated since the noncommutative Doob maximal inequality is now  involved. 
On the other hand, the proof of the duality theorem for $1<p<2$ is also slightly harder than that in the case $p=1$. 
This partly explains why we have decided to first consider the case $p=1$. 

Let us recall the definition of the spaces $L_p(\M;\ell_\infty), 1\leq p \leq \infty$. 
A sequence $(x_n)_{n\geq 1}$ in $L_p(\M)$ belongs to $L_p(\M;\ell_\infty)$ 
if  $(x_n)_{n\geq 1}$ admits a factorization $x_n=ay_nb$ with $a,b \in L_{2p}(\M)$ and $(y_n)_{n\geq 1}\in \ell_\infty(L_\infty(\M))$. 
The norm of $(x_n)_{n\geq 1}$ is then defined as
$$\|(x_n)_{n\geq 1}\|_{L_p(\M;\ell_\infty)}=\inf_{x_n=ay_nb}\|a\|_{2p}\sup_{n\geq 1}\|y_n\|_\infty\|b\|_{2p}.$$
One can check that $(L_p(\M;\ell_\infty),\|\;\|_{L_p(\M;\ell_\infty)})$ is a Banach space. 
It is proved in \cite{ju-doob} and \cite{jx-erg} that if $(x_n)_{n\geq 1}$ is a positive sequence in $L_p(\M;\ell_\infty)$, then
\begin{equation}\label{norme sup}
\|(x_n)_{n\geq 1}\|_{L_p(\M;\ell_\infty)}
=\sup\Big\{\sum_{n\geq 1}\tau(x_ny_n) : y_n\in L_{p'}^+(\M) \mbox{ and } \Big\|\sum_{n\geq 1}y_n\Big\|_{p'}\leq 1\Big\}.
\end{equation}
The norm of $L_p(\M;\ell_\infty)$ will be denoted by $\|\sup^+_n x_n\|_p$. 
We should warn the reader that $\|\sup^+_n x_n\|_p$ is just a notation since $\sup_n x_n$ does not take any sense in the noncommutative setting.

Now let $2<q \leq \infty$. We define the space
$$L_q^c\mo(\M)=\big\{a\in L_2(\M) : \|{\sup_{n\geq 1}}^+ \E_n(|a-\E_n(a)|^2)\|_{q/2} <\infty\big\}$$
equipped with the norm
$$\|a\|_{L_q^c\mo(\M)}=\max\Big(\|\E_1(a)\|_q \; , \; \Big(\|{\sup_{n\geq 1}}^+ \E_n(|a-\E_n(a)|^2)\|_{q/2}\Big)^{1/2}\Big).$$
Then $(L_q^c\mo(\M),\|\cdot\|_{L_q^c\mo(\M)})$ is a Banach space. 
Similarly, we set
$$L_q^r\mo(\M)=\{a : a^* \in L_q^c\mo(\M)\}$$
equipped with the norm 
$$\|a\|_{L_q^r\mo(\M)}=\|a^*\|_{L_q^c\mo(\M)}.$$
Note that if $q=\infty$, then $L_\infty^c\mo=\bmo^c$ and $L_\infty^r\mo=\bmo^r$. 
For convenience we denote $L_q^c\mo(\M), L_q^r\mo(\M)$ respectively by $L_q^c\mo, L_q^r\mo$.

The following duality holds:

\begin{theorem}\label{th:dual_hp}
Let $1\leq p <2$ and $q$ be the index conjugate to $p$. 
Then $(\h_p^c)^*=L_q^c\mo$ with equivalent norms.\\
More precisely,
\begin{enumerate}[\rm (i)]
\item Every $a \in L_q^c\mo$ defines a continuous linear functional on $\h_p^c$ by
\begin{equation}\label{crochet dualite hp}
\phi_a(x)=\tau(a^*x),\quad \forall x \in L_2(\M).
\end{equation}
\item Conversely, any $\phi \in (\h_p^c)^*$ is given as above by some $a \in L_q^c\mo$.\newline
Moreover
\begin{equation}\label{dualite hp} 
\lambda_p^{-1/2}\|a\|_{L_q^c\mo} \leq \| \phi_a\|_{(\h_p^c)^*} \leq \sqrt{2}\, \|a\|_{L_q^c\mo}
\end{equation}
where $\lambda_p >0$ is a constant depending only on $p$ and $\lambda_p=O(1)$ as $p \rightarrow 1$, 
$\lambda_p \leq C(2-p)^{-2}$ as $p\rightarrow 2$.
\end{enumerate}
Similarly, we have $(\h_p^r)^*=L_q^r\mo$, and $(\h_p)^*=L_q^c\mo\cap L_q^r\mo\cap \h_q^d$.
\end{theorem}

\pf 
We show only the duality equality $(\h_p^c)^*=L_q^c\mo$. 
To this end, we will adapt the proof of the corresponding duality result for $\H_p^c$ in \cite{jx-burk} for the first step. 
The second one is adapted from the proof of Theorem \ref{th:dual_h1}. \\
\\
\textbf{Step 1:} 
Let $a\in L_q^c\mo$ and $x$ be a finite $L_2$-martingale such that $\|x\|_{\h_p^c}\leq 1$. Let $s$ be the index conjugate to $\frac{q}{2}$. 
We consider
$$\tilde{s}_{c,n}(x)=\Big(\sum_{k=1}^n\E_{k-1}|dx_k|^2\Big)^{p/2s} 
\quad \mbox{ and } \quad 
\tilde{s}_{c}(x)=\Big(\sum_{k=1}^{\infty}\E_{k-1}|dx_k|^2\Big)^{p/2s}.$$
Then $\tilde{s}_{c,n}(x)\in L_s(\M_n)$ and by approximation we may assume that the $\tilde{s}_{c,n}(x)$'s are invertible. 
By the arguments in the proof of the duality between $\h_1^c$ and $\bmo^c$ in Theorem \ref{th:dual_h1} we have
$$\begin{array}{ccl}
|\phi_a(x)|&\leq& \Big[ \tau\Big(\sum_{n\geq 1} \tilde{s}_{c,n}(x)|da_n|^2\Big)\Big]^{1/2}
\Big[ \tau\Big(\sum_{n\geq1}\tilde{s}_{c,n}(x)^{-1/2}|dx_n|^2\tilde{s}_{c,n}(x)^{-1/2}\Big)\Big]^{1/2}\\
&=:& I\cdot II.
\end{array}$$
To estimate $I$ we set again
$$\left\{\begin{array}{l}
\theta_1=\tilde{s}_{c,1}(x)\\
\theta_n=\tilde{s}_{c,n}(x)-\tilde{s}_{c,n-1}(x), \quad \forall n\geq 2.
\end{array}\right.$$
Then $\theta_n \in L_s(\M_{n-1})$, $\theta_n \geq 0$ and 
$\tilde{s}_{c,n}(x)=\displaystyle \sum_{k=1}^n \theta_k$. 
Thus
$$\Big\| \sum_{k=1}^{\infty} \theta_k\Big\|_s=\|\tilde{s}_{c}(x)\|_s=\|x\|_{h_p^c}^{p/s}\leq 1.$$ 
By (\ref{norme sup}), we have
$$\begin{array}{ccl}
I^2&=& \displaystyle \sum_{k\geq1}\tau\Big(\theta_k\displaystyle \sum_{n\geq k}|da_n|^2\Big) \\
&=&\displaystyle \sum_{k\geq1}\tau\Big(\theta_k\E_{k-1}\Big(\displaystyle \sum_{n\geq k}|da_n|^2\Big)\Big) \\
&=&\displaystyle \sum_{k\geq 1} \tau(\theta_k\E_{k-1}(|a-a_{k-1}|^2))\\
&\leq &\|\displaystyle {\sup_{k \geq1}}^+\E_k(|a-a_{k}|^2)\|_{q/2} = \|a\|^2_{L_q^c\mo}.
\end{array}$$
To estimate the second term, let $\alpha=2/p \in (1,2]$ and notice that
$$1-\alpha=1-\frac{2}{p}=1-2+\frac{2}{q}=-\frac{1}{s}.$$
For fixed $n$, we define $y=\tilde{s}_{c,n-1}(x)^s$ and $z= \tilde{s}_{c,n}(x)^s$. Since $p/2\leq 1$, we have
$$y=\Big(\sum_{k=1}^{n-1}\E_{k-1}|dx_k|^2\Big)^{p/2} \leq \Big(\sum_{k=1}^n\E_{k-1}|dx_k|^2\Big)^{p/2}=z.$$
Note that 
$$z^{\frac{1-\alpha}{2}}=z^{-\frac{1}{2s}}=\tilde{s}_{c,n}(x)^{-\frac{1}{2}}.$$ 
Applying Lemma $4.1$ of \cite{ju-doob}, we find
$$\begin{array}{ccl}
\tau(\tilde{s}_{c,n}(x)^{-1/2}\E_{n-1}|dx_n|^2\tilde{s}_{c,n}(x)^{-1/2})&=&\tau(z^{\frac{1-\alpha}{2}}(z^\alpha-y^\alpha)z^{\frac{1-\alpha}{2}})\\
&\leq& 2\tau(z-y)\\
&=&2\tau(\tilde{s}_{c,n}(x)^s-\tilde{s}_{c,n-1}(x)^s).
\end{array}$$
Therefore
$$\begin{array}{ccl}
II^2&\leq &2\displaystyle\sum_{n\geq1}\tau[\tilde{s}_{c,n}(x)^s-\tilde{s}_{c,n-1}(x)^s]\\
&=&2\tau[(\tilde{s}_{c}(x)^s]\\
&=&2\tau\Big[\Big(\displaystyle\sum_{n\geq 1}\E_{n-1}|dx_n|^2\Big)^{p/2}\Big]\\
&=&2\|x\|_{\h_p^c}^p \leq 2.
\end{array}$$
Combining the precedent estimations we deduce that for any finite $L_2$-martingale $x$ 
$$|\phi_a(x)| \leq \sqrt{2} \, \|a\|_{L_q^c\mo}\|x\|_{\h_p^c}.$$
Thus $\phi_a$ extends to an element of $(\h_p^c)^*$ with norm $\leq \sqrt{2}\, \|a\|_{L_q^c\mo}$.\\
\\
\textbf{Step 2:}
Let $\phi \in (\h_p^c)^*$ such that $\|\phi\|_{(\h_p^c)^*} \leq 1$. 
As $L_2(\M)\subset h_p^c$, $\phi$ induces a continuous functional $\tilde{\phi}$ on $L_2(\M)$. 
Thus there exists $a \in L_2(\M)$ such that 
$$ \tilde{\phi}(x)=\tau(a^*x),\quad \forall x \in L_2(\M).$$
By the density of $L_2(\M)$ in $\h_p^c$ we have 
\begin{equation}\label{norme phi p}
\|\phi\|_{(\h_p^c)^*}=\sup_{x\in L_2(\M), \|x\|_{\h_p^c}\leq 1} |\tau(a^*x)|\leq 1.
\end{equation}
We want to estimate 
$$\|a\|_{L_q^c\mo}^2=
\max\Big(\|\E_1(a)\|_q^2 \; , \; \Big\|{\sup_{n\geq 1}}^+ \E_n\Big(\sum_{k>n}|da_k|^2\Big)\Big\|_{q/2}\Big).$$
Let $x\in L_p(\M_1), \|x\|_p\leq 1$ be such that $\|\E_1(a)\|_q=|\tau(a^*x)|$. 
Then by \eqref{norme phi p} we have 
$$\|\E_1(a)\|_q\leq \|x\|_{\h_p^c}=\|x\|_p\leq 1.$$
On the other hand for each $n\geq 1$ we set 
$$z_n=s_c(a)^2-s_{c,n}(a)^2=\sum_{k>n}\E_{k-1}|da_k|^2.$$
Then by (\ref{norme sup}) and the dual form of Junge's noncommutative Doob maximal inequality, 
we find (recalling that $s$ is the conjugate index of $q/2$) 
$$\begin{array}{cl}
&\|{\displaystyle\sup_{n\geq 1}}^+ \E_n(z_n)\|_{q/2}\\
=&\sup\Big\{\displaystyle\sum_{n\geq 1}\tau(\E_n(z_n)b_n):
b_n \in L_s^+(\M) \mbox{ and } \Big\|\displaystyle\sum_{n\geq 1} b_n\Big\|_s\leq 1\Big\}\\
\leq &\lambda_s \sup\Big\{\displaystyle\sum_{n\geq 1}\tau(\E_n(z_n)b_n):
b_n \in L_s^+(\M_n) \mbox{ and } \Big\|\displaystyle\sum_{n\geq 1} b_n\Big\|_s\leq 1\Big\}.
\end{array}$$
Note that $\lambda_s=O(1)$ as $s$ close to $1$, 
so $\lambda_s$ remains bounded as $q\rightarrow \infty$, i.e, as $p\rightarrow 1$. 
On the other hand, $\lambda_s \approx s^2$ as $s\rightarrow \infty$, i.e, as $p\rightarrow 2$.

Let $(b_n)_{n\geq 1}$ be a sequence in $L_s^+(\M_n)$ 
such that $\Big\|\displaystyle\sum_{n\geq 1} b_n\Big\|_s\leq 1$. 
Let $y$ be the martingale defined as follows 
$$dy_k=da_k\Big(\sum_{k>n}b_n\Big), \quad \forall k\geq 1.$$
By (\ref{norme phi p}) we have
$$\tau(a^*y) \leq \|y\|_{\h_p^c}.$$
Since $b_n\in L_s^+(\M_n)$ for any $n\geq 1$, we have
$$\begin{array}{ccl}
\tau(a^*y)&=&\tau\Big(\displaystyle\sum_{k\geq 1} da_k^*dy_k\Big)
=\tau\Big(\displaystyle\sum_{k\geq 1} |da_k|^2\Big(\displaystyle\sum_{k>n}b_n\Big)\Big)\\
&=&\displaystyle\sum_{n\geq 1}\displaystyle\sum_{k>n}\tau(|da_k|^2b_n)
=\displaystyle\sum_{n\geq 1}\displaystyle\sum_{k>n}\tau(\E_{k-1}|da_k|^2b_n)\\
&=&\displaystyle\sum_{n\geq 1}\tau(z_nb_n)\\
&=&\displaystyle\sum_{n\geq 1}\tau(\E_n(z_n)b_n).
\end{array}$$
On the other hand, by the definition of $y$ and the fact that $b_n\in L_1^+(\M_n)$, we find
$$\begin{array}{ccl}
s_c(y)^2&=&\displaystyle\sum_{k\geq 1} \E_{k-1}|dy_k|^2
=\displaystyle\sum_{k\geq 1} \E_{k-1}\Big|da_k \Big(\sum_{k>n}b_n\Big)\Big|^2\\
&=&\displaystyle\sum_{k\geq 1} \E_{k-1}\Big[\Big(\sum_{k>n}b_n\Big)|da_k|^2\Big(\sum_{k>n}b_n\Big)\Big]
=\displaystyle\sum_{k\geq 1}\displaystyle\sum_{n,m<k} b_n\E_{k-1}(|da_k|^2)b_m\\
&=&\displaystyle\sum_{n,m\geq 1} b_nz_{\max(n,m)}b_m.
\end{array}$$
We consider the tensor product $\mathcal{N}=\M\overline{\otimes} B(\ell_2)$, 
equipped with the trace $\tau \otimes \mathrm{tr}$, where $ \mathrm{tr}$ denote the usual trace on $B(\ell_2)$. 
Note that
$$\begin{array}{cl}
&\left[\begin{array}{ccc}
b_1^{1/2}& b_2^{1/2} & \hdots \\
0 & 0 &\hdots \\
\vdots & \vdots & 
\end{array}\right]
\Big[ b_n^{1/2}z_{\max(n,m)}b_m^{1/2}\Big]_{n,m\geq 1}
\left[\begin{array}{ccc}
b_1^{1/2}& 0 & \hdots \\
b_2^{1/2}& 0 &\hdots \\
\vdots & \vdots & 
\end{array}\right]\\
=&
\left[\begin{array}{ccc}
\displaystyle\sum_{n,m\geq 1} b_nz_{\max(n,m)}b_m & 0 & \hdots \\
0 & \hdots & \hdots \\
\vdots & &
\end{array}\right].
\end{array}$$
We claim that the matrix $Z=\Big[ b_n^{1/2}z_{\max(n,m)}b_m^{1/2}\Big]_{n,m\geq 1}$ is positive. 
Indeed, we suppose that $\M$ acts on the Hilbert space $H$ and we denote by $\langle \cdot, \cdot \rangle$ the associated scalar product. 
For $\xi=(\xi_n)_{n\geq 1} \in \ell_2(H)$, we have
$$\begin{array}{ccl}
\langle Z\xi, \xi \rangle_{\ell_2(H)} & = &
\displaystyle\sum_{n,m\geq 1}\langle \big(b_n^{1/2}z_{\max(n,m)}b_m^{1/2}\big)\xi_m,\xi_n\rangle\\
&=&\displaystyle\sum_{n,m\geq 1}\langle z_{\max(n,m)}\big(b_m^{1/2}\xi_m\big),b_n^{1/2}\xi_n\rangle,
\end{array}$$
where the last equality comes from the positivity of the $b_n$'s. 
Then the definition of $z_n$ gives
$$\begin{array}{ccl}
\langle Z\xi, \xi \rangle_{\ell_2(H)} & = &
\displaystyle\sum_{n,m\geq 1}
\langle \Big(\displaystyle\sum_{k>\max(n,m)}\E_{k-1}|da_k|^2\Big)\big(b_m^{1/2}\xi_m\big),b_n^{1/2}\xi_n\rangle \\
&=& \displaystyle\sum_{k\geq 1}
\langle \E_{k-1}|da_k|^2\Big(\displaystyle\sum_{m<k}b_m^{1/2}\xi_m\Big),
\Big(\displaystyle\sum_{n<k}b_n^{1/2}\xi_n\Big)\rangle. 
\end{array}$$
The positivity of the conditional expectation implies that each term of the latter sum is non-negative. 
Thus, we obtain 
$$\langle Z\xi, \xi \rangle_{\ell_2(H)} \geq 0, \quad \forall \xi \in \ell_2(H),$$
which proves our claim. 
Hence 
$$\|Z\|_{L_{1}(\mathcal{N})}=\tau \otimes \mathrm{tr}(Z)=\sum_{n\geq 1}\tau\Big(b_n^{1/2}z_nb_n^{1/2}\Big).$$
Since $\frac{2}{p}=\frac{1}{2s}+1+\frac{1}{2s}$, by the H\"{o}lder inequality we have 
$$\begin{array}{cl}
&\Big\|\displaystyle\sum_{n,m\geq 1} b_nz_{\max(n,m)}b_m \Big\|_{p/2} 
=\Big\|\Big(\displaystyle\sum_{n,m\geq 1} b_nz_{\max(n,m)}b_m \Big)\otimes e_{1,1}\Big\|_{L_{p/2}(\mathcal{N})}\\
\leq &
\Big\|\left[\begin{array}{ccl}
b_1^{1/2}& b_2^{1/2} & \hdots \\
0 & 0 &\hdots \\
\vdots & \vdots & 
\end{array}\right]\Big\|_{L_{2s}(\mathcal{N})}
\Big\|\Big[ b_n^{1/2}z_{\max(n,m)}b_m^{1/2}\Big]_{n,m\geq 1}\Big\|_{L_{1}(\mathcal{N})}
\Big\|\left[\begin{array}{ccc}
b_1^{1/2}& 0 & \hdots \\
b_2^{1/2}& 0 &\hdots \\
\vdots & \vdots & 
\end{array}\right]\Big\|_{L_{2s}(\mathcal{N})}\\
=& \Big\|\displaystyle\sum_{n\geq 1} b_n\Big\|_s^{1/2} 
\Big[\displaystyle\sum_{n\geq 1}\tau\Big(b_n^{1/2}z_nb_n^{1/2}\Big)\Big]
\Big\|\displaystyle\sum_{n\geq 1} b_n\Big\|_s^{1/2}\\
=&\Big[\displaystyle\sum_{n\geq 1}\tau(\E_n(z_n)b_n)\Big] \Big\|\displaystyle\sum_{n\geq 1} b_n\Big\|_s.
\end{array}$$
Thus 
$$\|y\|_{\h_p^c}\leq \Big[\displaystyle\sum_{n\geq 1}\tau(\E_n(z_n)b_n)\Big]^{1/2} \Big\|\displaystyle\sum_{n\geq 1} b_n\Big\|_s^{1/2}.$$
Combining the preceding inequalities, we deduce
$$\sum_{n\geq 1}\tau(\E_n(z_n)b_n)\leq \Big\|\displaystyle\sum_{n\geq 1} b_n\Big\|_s\leq 1.$$
Therefore $a\in L_q^c\mo$ and $\|a\|_{L_q^c\mo}\leq \sqrt{\lambda_s}$. 
This ends the proof of the duality $(\h_p^c)^*=L_q^c\mo$. 
\cqd

\vspace{0.5cm}

\begin{rk}\label{rk:duality-cst}
Junge and Mei obtain in \cite{jm-riesz} the following inequality
$$\lambda_p^{-1}\|a\|_{L_q^c\mo} \leq \| \phi_a\|_{(\h_p^c)^*} \leq \sqrt{2}\, \|a\|_{L_q^c\mo}$$
where $\lambda_p$ is the constant in (\ref{dualite hp}).
Note that our lower estimate is the square root of theirs, and yields a better estimation as $p\rightarrow 2$.
\end{rk}

The dual space of $\H_p$ for $1\leq p <2$ is described in \cite{jx-burk} 
as the space $L_q\MO$ (where $q$ is the index conjugate of $p$) defined as follows. 
Let $2<q\leq \infty$, we set 
$$L_q^c\MO(\M)=\big\{a\in L_2(\M) : \|{\sup_{n\geq 1}}^+ \E_n(|a-\E_{n-1}(a)|^2)\|_{q/2} <\infty\big\},$$
equipped with the norm
$$\|a\|_{L_q^c\MO(\M)}=\Big(\|{\sup_{n\geq 1}}^+ \E_n(|a-\E_{n-1}(a)|^2)\|_{q/2}\Big)^{1/2}.$$
Similarly, we define
$$L_q^r\MO(\M)=\{a : a^* \in L_q^c\MO(\M)\},$$
equipped with the norm 
$$\|a\|_{L_q^r\MO(\M)}=\|a^*\|_{L_q^c\MO(\M)}.$$
Finally, we set 
$$L_q\MO(\M)=L_q^c\MO(\M) \cap L_q^r\MO(\M),$$
equipped with the intersection norm
$$\|x\|_{L_q\MO(\M)}=\max\big(\|x\|_{L_q^c\MO(\M)} \; , \; \|x\|_{L_q^r\MO(\M)}\big).$$
Note that if $q=\infty$, these spaces coincide with the $\BMO$ spaces. 
For convenience we denote $L_q^c\MO(\M), L_q^r\MO(\M), L_q\MO(\M)$ 
respectively by $L_q^c\MO, L_q^r\MO, L_q\MO$.\\
Theorem $4.1$ of \cite{jx-burk} establishes the duality $(\H_p^c)^*=L_q^c\MO$.  
Moreover, for any $a\in L_q^c\MO$,
$$\lambda_p^{-1}\|a\|_{L_q^c\MO} \leq \| \phi_a\|_{(\H_p^c)^*} \leq \sqrt{2}\, \|a\|_{L_q^c\MO}$$
where $\lambda_p$ is the constant in (\ref{dualite hp}).

\begin{rk} 
The method used in the second step of the previous proof can be adapted to the duality between 
$\H_p^c$ and $L_q^c\MO$, for $1< p <2$. 
This yields a better estimate of the constant $\lambda_p$ given in \cite{jx-burk}. 
More precisely, we obtain by this way a constant of order $(2-p)^{-1}$ as $p\rightarrow 2$ 
instead of $(2-p)^{-2}$. 

Indeed, let $\phi \in (\H_p^c)^*$ such that $\|\phi\|_{(\H_p^c)^*} \leq 1$. 
There exists $a \in L_2(\M)$ such that 
$$ \phi(x)=\tau(a^*x),\quad \forall x \in L_2(\M).$$
By the density of $L_2(\M)$ in $\H_p^c$ we have 
$$\|\phi\|_{(\H_p^c)^*}=\sup_{x\in L_2(\M), \|x\|_{\H_p^c}\leq 1} |\tau(a^*x)|\leq 1.$$
In this case, we want to estimate 
$$\|a\|_{L_q^c\MO}^2=\big\|{\sup_{n\geq 1}}^+ \E_n(|a-\E_{n-1}(a)|^2)\big\|_{q/2}
=\Big\|{\sup_{n\geq 1}}^+ \E_n\Big(\sum_{k\geq n}|da_k|^2\Big)\Big\|_{q/2}.$$
The triangular inequality in $L_{q/2}(\M;\ell_\infty)$ 
allows us to separate the estimation into two parts as follows
$$\|a\|_{L_q^c\MO}^2 \approx 
\Big\|{\sup_{n\geq 1}}^+ \E_n\Big(\sum_{k>n}|da_k|^2\Big)\Big\|_{q/2} 
+ \big\|{\sup_{n\geq 1}}^+ |da_n|^2\big\|_{q/2} 
=: I +II.$$
We adapt the second step of the preceding proof by 
setting $z_n=\displaystyle\sum_{k>n}|da_k|^2$ for each $n\geq 1$. 
It yields the following estimation of the first term 
$$I\leq \lambda_s$$
where $s$ is the index conjugate to $\frac{q}{2}$.

To estimate the diagonal term $II$, 
let $(b_n)_{n\geq 1}$ be a sequence in $L_s^+(\M_n)$ 
such that \\
$\Big\|\displaystyle\sum_{n\geq 1} b_n\Big\|_s\leq 1$. 
Let $y$ be the martingale defined as follows 
$$dy_k=da_kb_k-\E_{k-1}(da_kb_k), \quad \forall k\geq 1.$$
We have
$$\tau(a^*y) \leq \|y\|_{\H_p^c}.$$
Since $(da_n)_{n\geq 1}$ is a martingale difference sequence, we have
$$\begin{array}{ccl}
\tau(a^*y)&=&\displaystyle\sum_{n\geq 1} \tau(|da_n|^2b_n)-\tau(da_n^*\E_{n-1}(da_nb_n))\\
&=&\displaystyle\sum_{n\geq 1} \tau(|da_n|^2b_n)-\tau(\E_{n-1}(da_n^*)da_nb_n)\\
&=&\displaystyle\sum_{n\geq 1} \tau(|da_n|^2b_n).
\end{array}$$
On the other hand, the triangular inequality in $L_p(\M;\ell_2^c)$ yields
$$\|y\|_{\H_p^c}=\|(dy_n)_{n\geq 1}\|_{L_p(\M;\ell_2^c)}
\leq
\big\|\big(da_nb_n\big)_{n\geq 1}\big\|_{L_p(\M;\ell_2^c)}
+\big\|\big(\E_{n-1}(da_nb_n)\big)_{n\geq 1}\big\|_{L_p(\M;\ell_2^c)}$$
The noncommutative Stein inequality implies
$$\big\|\big(\E_{n-1}(da_nb_n)\big)_{n\geq 1}\big\|_{L_p(\M;\ell_2^c)}
\leq\gamma_p \big\|\big(da_nb_n\big)_{n\geq 1}\big\|_{L_p(\M;\ell_2^c)}$$
with $\gamma_p\leq C\frac{p^2}{p-1}$ (see \cite{jx-burk}). 
Then 
$$\|y\|_{\H_p^c}\leq
(1+\gamma_p)\big\|\big(da_nb_n\big)_{n\geq 1}\big\|_{L_p(\M;\ell_2^c)}.$$
As before, by the H\"{o}lder inequality, we find
$$\begin{array}{cl}
&\big\|\big(da_nb_n\big)_{n\geq 1}\big\|_{L_p(\M;\ell_2^c)}^2
=\Big\|\displaystyle\sum_{n\geq 1} b_n|da_n|^2b_n \Big\|_{p/2} \\
\leq &
\Big\|\displaystyle\sum_{n\geq 1} b_n\Big\|_s^{1/2} 
\Big[\displaystyle\sum_{n\geq 1}\tau\Big(b_n^{1/2}|da_n|^2b_n^{1/2}\Big)\Big]
\Big\|\displaystyle\sum_{n\geq 1} b_n\Big\|_s^{1/2}\\
=&\Big[\displaystyle\sum_{n\geq 1}\tau(|da_n|^2b_n)\Big] \Big\|\displaystyle\sum_{n\geq 1} b_n\Big\|_s.
\end{array}$$
Hence
$$\|y\|_{\H_p^c}\leq
(1+\gamma_p)
\Big\|\sum_{n\geq 1} b_n\Big\|_s^{1/2}\Big(\sum_{n\geq 1}\tau(|da_n|^2b_n)\Big)^{1/2}.$$
Combining the preceding inequalities, we deduce
$$\sum_{n\geq 1}\tau(|da_n|^2b_n)
\leq C(p-1)^{-2}\Big\|\sum_{n\geq 1} b_n\Big\|_s.$$
Then
$$II=\big\|{\sup_{n\geq 1}}^+ |da_n|^2\big\|_{q/2} \leq C(p-1)^{-2}\lambda_s.$$
Finally, we obtain 
$$\|a\|_{L_q^c\MO}\leq \lambda_s^{1/2}(1+C(p-1)^{-2})^{1/2}.$$
Since $\lambda_s\approx s^2$ as $s\rightarrow \infty$, i.e, as $p\rightarrow 2$, 
we have the announced estimation $\lambda_s^{1/2}(1+C(p-1)^{-2})^{1/2}\approx (2-p)^{-1}$ as $p\rightarrow 2$.
\end{rk}

\vspace{0.5cm}

For $1<p<\infty$, the noncommutative Burkholder-Gundy inequalities
of \cite{px-BG} 
and the noncommutative Burkholder inequalities of \cite{jx-burk} 
state respectively that $\H_p(\M)=L_p(\M)$ and $\h_p(\M)=L_p(\M)$ (with equivalent norms). 
Combining these results we obtain the equivalence of the norms $\H_p$ and $\h_p$. 
This is stated in Proposition $6.2$ of \cite{ran-cond}. 
Here Theorem \ref{th:dual_hp} allows us to compare the dual spaces of $\H_p$ and $\h_p$ for $1\leq p <2$. 
This dual approach gives another way to compare the spaces $\H_p$ and $\h_p$ for $1\leq p <2$, 
which improve the estimation of the constant $\kappa_p$ below for $1<p<2$. 
Indeed, Randrianantoananina obtained $\kappa_p=O((p-1)^{-1})$ as $p\rightarrow 1$ 
and the following statement gives that $\kappa_p$ remains bounded as $p\rightarrow 1$.  
For completeness, we also include Randrianantoanina's estimates.
 
\begin{theorem}\label{th:Hp=hp}
Let $1<p<\infty$. There exist two constants $\kappa_p>0$ and $\nu_p>0$ (depending only on $p$) 
such that for any finite $L_p$-martingale $x$,
$$\kappa_p^{-1}\|x\|_{\h_p}\leq  \|x\|_{\H_p}\leq \nu_p\|x\|_{\h_p}.$$
Moreover 
\begin{enumerate}[\rm (i)]
\item $\kappa_p \approx 1$ as $p\rightarrow 1$; 
\item $\kappa_p \leq Cp$ for $2\leq p <\infty$; 
\item $\nu_p \approx 1$ as $p\rightarrow 1$; 
\item $\nu_p \leq C\sqrt{p}$ for $2\leq p <\infty$.
\end{enumerate}
\end{theorem}

\vspace{0.5cm}

\pf 
Randrianantoanina stated the estimations $(ii), (iii), (iv)$ in  \cite{ran-cond} without giving the proof. 
For the sake of completness we give the proof of these three estimations.

\textbf{(i)}
Here we adopt a dual approch. 
Let $1 < p <2$ and $q$ the index conjugate to $p$. 
Let $a \in (\h_p)^*$. Then the triangular inequality in 
$L_{q/2}(\M;\ell_\infty)$ gives
$$\begin{array}{ccl}
\|a\|_{L_q^c\MO}^2&=&\Big\|\displaystyle{\sup_{n\geq 1}}^+ \E_n\Big(\displaystyle\sum_{k\geq n} |da_k|^2\Big)\Big\|_{q/2}\\
&\leq&\Big\|\displaystyle{\sup_{n\geq 1}}^+ \E_n\Big(\displaystyle\sum_{k> n} |da_k|^2\Big)\Big\|_{q/2}
+\|\displaystyle{\sup_{n\geq 1}}^+ \E_n|da_n|^2\|_{q/2}\\
&\leq& \|a\|_{L_q^c\mo}^2+\|(\E_n|da_n|^2)_n\|_{L_{q/2}(\M,\ell_\infty)}
\end{array}$$
But for $1\leq p \leq \infty$ we have the following contractive inclusion
$$\ell_p(L_p(\M))\subset L_p(\M;\ell_\infty).$$
Therefore
$$\begin{array}{ccl}
\|(\E_n|da_n|^2)_n\|_{L_{q/2}(\M;\ell_\infty)} &\leq &\|(\E_n|da_n|^2)_n\|_{\ell_{q/2}(L_{q/2})}\\
&\leq& \Big(\displaystyle\sum_{n\geq 1} \|da_n\|_q^q\Big)^{2/q}
=\|a\|_{\h_q^d}^2.
\end{array}$$
Then
$$\begin{array}{ccl}
2^{-1/2}\|a\|_{(\H_p)^*}&\leq &\|a\|_{L_q\MO}\\
&\leq& \sqrt{2} \max(\|a\|_{L_q^c\mo},\|a\|_{L_q^r\mo},\|a\|_{h_q^d}) \\
&\leq &\sqrt{2\lambda_p}\;\|a\|_{(h_p)^*}
\end{array}$$
with $\lambda_p=O(1)$ as $p\rightarrow 1$, hence $\kappa_p \approx 1$ as $p\rightarrow 1$.

\textbf{(ii)} The dual version of the noncommutative Doob inequality in \cite{ju-doob} 
gives that for $1\leq p < \infty$ and for all finite sequences $(a_n)$ of positive elements in $L_p(\M)$ :
$$\Big\|\sum_{n\geq 1} \E_{n-1}(a_n)\Big\|_p\leq c_p\Big\|\sum_{n\geq 1} a_n\Big\|_p$$
with $c_p \approx p^2$ as $p\rightarrow +\infty$. 
Applying this to $a_n=|dx_n|^2$ and $p/2$ we get
$$\begin{array}{ccl}
\|x\|_{h_p^c} &=    &\Big\|\Big(\displaystyle\sum_{n\geq 1}\E_{n-1}|dx_n|^2\Big)^{1/2}\Big\|_p 
                     = \Big\|\displaystyle\sum_{n\geq 1}\E_{n-1}|dx_n|^2\Big\|_{p/2}^{1/2}\\
              &\leq &\sqrt{c_{p/2}}\, \Big\|\displaystyle\sum_{n\geq 1}|dx_n|^2\Big\|_{p/2}^{1/2}=\sqrt{c_{p/2}}\,\|x\|_{\H_p^c}.
\end{array}$$
Passing to adjoints we have $\|x\|_{h_p^r}\leq \sqrt{c_{p/2}}\,\|x\|_{\H_p^r}$ with $\sqrt{c_{p/2}}\approx p$ as $p\rightarrow \infty$.\\
On the other hand, we have for $2\leq p \leq \infty$ and for any finite sequence $(a_n)$ in $L_p(\M)$
$$\Big(\sum_{n\geq 1}\|a_n\|_p^p\Big)^{1/p}\leq \Big\|\Big(\sum_{n\geq1}|a_n|^2\Big)^{1/2}\Big\|_p.$$
Indeed, this is trivially true for $p=2$ and $p=\infty$. 
Then complex interpolation yields the intermediate case $2<p<\infty$. \\
It thus follows that $\|dx\|_{\ell_p(L_p)}\leq \|x\|_{\H_p^c}$.\\
Thus $\kappa_p \leq Cp$ for $2\leq p <\infty$.

\textbf{(iii)} Adapting the discussion following Theorem \ref{th:H1=h1} to the case $0<p<1$, we obtain this estimate.

\textbf{(iv)} Suppose $2<p<\infty$ and $\|x\|_{\h_p} \leq 1$. We write
$$|dx_n|^2=\E_{n-1}|dx_n|^2+(|dx_n|^2-\E_{n-1}|dx_n|^2)=:\E_{n-1}|dx_n|^2+dy_n.$$
The noncommutative Burkholder inequality implies
$$\begin{array}{ccl}
\Big\|\displaystyle\sum_{n\geq 1}dy_n\Big\|_{p/2}&=&\|y\|_{p/2}\\
&\leq &\eta_{p/2} \Big[\Big(\displaystyle\sum_{n\geq 1}\|dy_n\|_{p/2}^{p/2}\Big)^{2/p}+\Big\|\Big(\displaystyle\sum_{n\geq 1}\E_{n-1}|dy_n|^2\Big)^{1/2}\Big\|_{p/2}\Big]
=: \eta_{p/2}(I+II)
\end{array}$$
with $ \eta_{p/2}\leq Cp$ as $p\rightarrow \infty$ from the proof of Theorem $4.1$ of \cite{ran-cond}. 
In order to estimate $I$ we use the triangular inequality in $\ell_{p/2}(L_{p/2})$ and contractivity of the conditional expectations:
$$I=\|dy\|_{\ell_{p/2}(L_{p/2})}\leq 2 \Big(\sum_{n\geq 1}\||dx_n|^2\|_{p/2}^{p/2}\Big)^{2/p}=2\Big(\sum_{n\geq 1}\|dx_n\|_{p}^{p}\Big)^{2/p}\leq 2.$$
As for the second term $II$ we note that
$$\E_{n-1}|dy_n|^2
=\E_{n-1}|dx_n|^4-(\E_{n-1}|dx_n|^2)^2
\leq \E_{n-1}|dx_n|^4.$$
Then Lemma $5.2$ of \cite{jx-burk} gives the following estimation
$$\begin{array}{ccl}
II&\leq & \Big\|\Big(\displaystyle\sum_{n\geq 1}\E_{n-1}|dx_n|^4\Big)^{1/2}\Big\|_{p/2}
= \Big\|\displaystyle\sum_{n\geq 1}\E_{n-1}|dx_n|^4\Big\|_{p/4}^{1/2}\\
&\leq & \Big(\Big\|\displaystyle\sum_{n\geq 1}\E_{n-1}|dx_n|^2\Big\|_{p/2}^{(p-4)/(p-2)} 
\Big(\displaystyle\sum_{n\geq 1}\|dx_n\|_{p}^{p}\Big)^{2/(p-2)}\Big)^{1/2}\\
&\leq & 1.
\end{array}$$
Combining the preceding inequalities we obtain
$$\|x\|_{\H_p^c}^2\leq 1+3\eta_{p/2}\leq Cp \mbox{ as } p\rightarrow \infty.$$
\cqd

\vspace{0.5cm}

\begin{rk}
At the time of this writing, we do not know if the orders of growth of $\kappa_p$ and $\nu_p$ for $2<p<\infty$ are optimal.
\end{rk}

\addcontentsline{toc}{chapter}{Bibliography}
\nocite{ran-mtrans}

\end{document}